\documentclass[a4paper,reqno,12pt]{amsart}

\usepackage{amsthm,amsfonts,amsxtra,amssymb,amscd}

\usepackage[all]{xy}

\theoremstyle{plain}       
\newtheorem{lem}{Lemma}
\newtheorem{teo}[lem]{Theorem}
\newtheorem{prp}[lem]{Proposition}
\newtheorem{cor}[lem]{Corollary}

\theoremstyle{definition}
\newtheorem{dfn}[lem]{Definition}
\newtheorem{oss}[lem]{Remark}

\theoremstyle{remark}
\newtheorem{ntz}[lem]{Notation}

\newcommand{\mC}  {{\mathbb C}}  

\newcommand{\mN}  {{\mathbb N}}

\newcommand{\calN}  {\text{$\mathcal {N}$}}

\newcommand{\calF}  {\text{$\mathcal {F}$}}

\newcommand{\calS}  {\text{$\mathcal {S}$}}
\newcommand{\calP}  {\text{$\mathcal {P}$}}

\newcommand{\calR}  {\text{$\mathcal {R}$}}
\newcommand{\calI}  {\text{$\mathcal {I}$}}

\newcommand{\grg}  {\text{$\gamma$}}    
\newcommand{\grd}  {\text{$\delta$}}

\newcommand{\grO}  {\text{$\Omega$}}
\newcommand{\gra}  {\text{$\alpha$}}
\newcommand{\grl}  {\text{$\lambda$}}
\newcommand{\grL}  {\text{$\Lambda$}}

\newcommand{\grs}  {\text{$\sigma$}}
\newcommand{\grb}  {\text{$\beta$}}

\newcommand{\grf}  {\text{$\varphi$}}

\newcommand{\jp}  {j^{\prime}}
\newcommand{\hp}  {h^{\prime}}

\newcommand{\vpri}  {v^{\prime}}

\newcommand{\trasp}[1]   {\sideset{^{t}}{}{\operatorname{#1}} }

\newcommand{\tc}         {\, : \,}

\newcommand{\degree}     {\operatorname{deg}}    
\newcommand{\grad}       {\operatorname{grad}}    
\newcommand{\Tr}         {\operatorname{Tr}}    
    
\renewcommand{\Im}         {\operatorname{Im}}    
\newcommand{\then}     {\Rightarrow}

\newcommand{\mfor}       {\text{ for }}
\newcommand{\mand}       {\text{ and }}
\newcommand{\mif}        {\text{ if }}

\newcommand{\Hom}        {\text{Hom}}

\newcommand{\wt}[1]      { {\widetilde {#1} } }
\newcommand{\wbar}[1]    { {\overline  {#1} } }

\newcommand{\ra}         {\rightarrow}
\newcommand{\lra}        {\longrightarrow}

\newcommand{\isocan} {\simeq}
\newcommand{\vuoto}  {\varnothing}

\newcommand{\slice}     {\mbox{${\mathcal{S}} _{a,x}$} }
\newcommand{\sopraslice} {\mbox{${\widetilde {\mathcal{S}} _{a,x}}$} }

\begin{document}

\title{Quiver varieties of type $\mathsf{A}$}
\author{Andrea Maffei}
\begin{abstract}  
We prove a conjecture of Nakajima describing the  relation between the geometry of
quiver varieties of type $\mathsf{A}$ 
and the geometry of the partial flags varieties 
and of the nilpotent 
variety.   
\end{abstract}

\maketitle

\bigskip
The kind of quiver varieties we are interested in, have been introduced by Nakajima as a generalization
of the description of the moduli space of anti-self-dual connections 
on ALE spaces constructed by Kroneheimer and Nakajima \cite{KrNa}. They result to have a rich and
interesting geometry 
and they were  used by Nakajima to give a geometric construction of the representations 
of Kac-Moody algebras \cite{Na1}, \cite{Na2}. A similar construction had already been done in the case of $sl_n$ by Ginzburg 
\cite{Ginzburg} using partial flags varieties.  A precise conjecture of Nakajima
(\cite{Na1} or theorem \ref{congettura} below) describes the relation between the two kind of varieties.

I want to thank Corrado De Concini who explained me quiver varieties and
pointed out to me this problem and Hiraku Nakajima who pointed out an error in
a previous version of this paper and the solution to it.

\section{Nakajima's conjecture}

We recall some definition and fix some notation on quiver varieties of type 
$\mathsf{A}_{n-1}$ and on partial flags varieties.

\subsection{Quiver varieties of type $\mathsf{A}_{n-1}$}
\label{convenzioniquiverA}

Let $C= 2I-A$ 
be the Cartan matrix of type $\mathsf{A}_{n-1}$ and consider the ``double''
graph of type $\mathsf{A}_{n-1}$. We will call the vertices  and the arrows 
of this graph according to the following diagram:
\begin{equation}\label{diag:tipoA}
   \xymatrix{
	1  \ar@/_/[r]_{ar_1} &  
	2  \ar@/_/[r]_{ar_2} \ar@/_/[l]_{\overline{ ar} _1} &
	{\cdots}\ar@/_/[l]_{\overline{ ar} _2} &
	n-2  \ar@/_/[r]_{ar_{n-2}}  &
	n-1  \ar@/_/[l]_{\overline{ ar} _{n-2}}}.
\end{equation}
In particular $I = \{1,\dots,n-1\}$ is the set of vertices,
$\Omega = \{ar_1,\dots,ar_{n-2}\}$, $\overline{\Omega}= \{\overline{ar} _1 , 
\dots , \overline{ ar} _{n-2} \}$  and $H=\Omega \cup \overline{\Omega}$ is the
set of arrows. We observe that given $i,j \in I $ we have that 
$a_{ij} = card \{h \in \Omega$ joining $i$ and $j\}$. Finally if $h$ is an 
arrow we call $h_0$ its source of $h$ and $h_1$ its target.
 
\begin{ntz} \label{notazioni1}
In this paper  
$v=(v_{1},\dots,v_{n-1})$ and $d=(d_{1},\dots,d_{n-1})$ 
will be vectors 
of integers. In the case that they are vectors of nonnegative integers
$V_i$ and $D_i$  will be vector spaces of dimension $v_i$ and 
$d_i$ for $i=1,\dots,n-1$. We need also to define an action of the Weyl group
$S_n$ of the Dynkin diagram of type $A_{n-1}$ on $d$, $v$. Let $\gra_1,\dots,
\gra_{n-1}$ (resp. $\omega_1,\dots,\omega_{n-1}$ ) be simple roots 
(resp. fundamental weights) for this diagram numbered according to \eqref{diag:tipoA}.

We identify $n-1$-tuples of integers $v=(v_1,\dots,v_{n-1})$ with elements of
the root lattice through:
$$
v \longmapsto \sum_{i=1}^{n-1} v_i \gra_i
$$
and $n-1$-tuples of integers $d=(d_1,\dots,d_{n-1})$ with elements of the weight 
lattice through:
$$
d \longmapsto \sum_{i=1}^{n-1} d_i \omega_i
$$
We will use $v$ (resp. $d$) to denote 
both an $n-1$-tuples of integers and an element of the root (resp. weight) 
lattice.
\end{ntz}

If $d,v$ are vectors of nonnegative integers
we introduce now the vector space of 
``double free representation'' of this graph
``with respect'' the vector spaces $D_i,V_i$:
\begin{align*}
	S_{\grO}(d,v) &= \bigoplus _{i \in I} \Hom(D_i,V_i) \oplus \bigoplus 
			_{h \in \grO} \Hom(V_{h_0},V_{h_1}),\\
	S_{\wbar \grO}(d,v) &=\bigoplus _{i \in I} \Hom(V_i,D_i) \oplus \bigoplus 
			_{h \in \wbar \grO} \Hom(V_{h_0},V_{h_1}), \\
	S(d,v)&= S_{\grO}(d,v) \oplus S_{\wbar \grO}(d,v) .
\end{align*}
More often, when it will not be ambiguous we will write $S$ instead of 
$S(d,v)$.

\medskip

\begin{ntz} \label{notazioni2}
Following Nakajima if $h \in H$ and $i\in I$ 
we will use $B_h$ to denote an element 
of $\Hom(V_{h_0}, V_{h_1})$ and 
$\gamma_i$, (resp. $\grd_i$) to denote elements of  $Hom (D_{i},V_{i})$ 
(resp. $Hom (V_{i},D_{i})$).

In the case of the graph of type $\mathsf{A}$ the notation turn to be simpler
if we use a different notation for the maps $B_h$ according to  
$h\in \Omega$ or not. So we will use also 
 $A_i$ to denote an element of $Hom (V_i,V_{i+1})$ and
 $B_i$ to denote an element of $Hom (V_{i+1},V_{i})$.
Finally $A= (A_1, \dots , A_{n-2})$,  $B= (B_1, \dots , B_{n-2})$,
$\gamma= (\gamma_1, \dots , \gamma_{n-1})$ and 
$\grd= (\grd_1, \dots , \grd_{n-1})$. 
Our conventions can be summarized in the 
following diagram:
\begin{equation*}
   \xymatrix{
	D_1 \ar[d]^{\grg_1}  &  D_2 \ar[d]^{\grg_2}  &  &  
	D_{n-2} \ar[d]^{\grg_{n-2}} & D_{n-1}\ar[d]^{\grg_{n-1}}\\
	V_1 \ar[d]^{\grd_1} \ar@/_/[r]_{A_1} &  
	V_2 \ar[d]^{\grd_2} \ar@/_/[r]_{A_2} \ar@/_/[l]_{B_1} &
	{\cdots}\ar@/_/[l]_{B_2} &
	V_{n-2} \ar[d]^{\grd_{n-2}} \ar@/_/[r]_{A_{n-2}}  &
	V_{n-1} \ar[d]^{\grg_{n-1}} \ar@/_/[l]_{B_{n-2}}\\
	D_1  &	D_2  &	&  D_{n-2}  &  D_{n-1}			
   }
\end{equation*}
\end{ntz}

There is a natural action of
the groups   $GL(V) = \prod  GL(V_i)$ and $GL(D) = \prod  GL(D_i)$ 
on $S$: if $g=(g_i) \in GL(V)$ we define 
$$g( (A_i),(B_i),(\grg_i),(\grd_i))=( (g_{i+1} A_i g_i^{-1})
(g_{i} B_i g_{i+1}^{-1}), ( g_i \grg_i), (\grd_i g_i^{-1})).$$
The action of $GL(D)$ is described in a similar way.

We will use also the 
following notations: $\grg_{j \ra i} = B_{i} \dots B_{j-1} \grg_j$ 
and $\grd _{j \ra i} = \grd_i A_{i-1} 
\dots A_{j}$.

\begin{dfn}[Nakajima , \cite{Na1}, \cite{Na2}] \label{definizione quiver}
An element $(A,B,\gamma,\delta)$ of $S$ is said to satisfy the 
ADHM equation or $admissible$ if it satisfies the following 
relations
\begin{align*}
   B_1 A_1 &= \grg _1 \grd_1, \\
   B_i A_i &= A_{i-1} B_{i-1} + \grg_i \grd_i \quad 
				\mfor  2 \leq i \leq n-2, \\
   0	   &=  A_{n-2} B_{n-2} + \grg_{n-1} \grd_{n-1}.
\end{align*}
We call $\grL(d,v)$ the set of all admissible data.
An admissible data is said to be $stable$ if each collection $U=(U_1,\dots,U_{n-1})$ of subspaces of 
$V$ containing $\Im \grg$ and invariant by the action 
of $A$ and $B$ (that is $A_i(U_i)\subset U_{i+1}$ and 
$B_i(U_{i+1})\subset U_i$) must be equal to $V$. We call 
$\grL^+(d,v)$ the set of stable admissible data.

Observe that $\grL$ and $\grL^+$ are invariant by the 
action of $GL(V)$, so following Nakajima
we can define the following two types of quiver varieties:
\begin{align*}
  M^0(d,v) &= \grL(v,w) /\!/ GL(V) \mbox{ and  $p_0$ is the projection 
from  $\grL$  to $M^0$},\\
  M(d,v) &= \grL^+(v,w) /\!/ GL(V)  \mbox{ and  $p$ is the projection 
from  $\grL^+$ to $M$}.
\end{align*}
The construction come equipped with an  action of the group $GL(W)$ on the two 
quiver varieties 
and with an equivariant   projective map $\pi : M(d,v) \longrightarrow M^0(d,v)$. 
We observe that the data
$0=(0,0,0,0)$ is always admissible so we can define 
$\Gamma(d,v)$ to be the fiber of the map $\pi$ in zero.
This variety is called the Lagrangian subvariety. We define $M^1(d,v)$ 
to be the image of $\pi$,
which is closed since 
$\pi$ is projective, with 
the reduced structure.

For our purpose it will turn to be useful to extend the definition
of $M(d,v)$ and $M^0(d,v)$ also to the case of arbitrary 
$n-1$-tuples of integers: $M(d,v)=M^0(d,v)=\vuoto$
if there exists $i$ such that $v_i < 0 $ or $d_i <0$.
\end{dfn}

\subsection{Admissible polynomials and the coordinate ring of $M^0$}

We have defined $M^0$ as the affine variety $\grL /\!/G$.
A set of generator of its coordinate ring has been given by Lusztig in
\cite{Lu:Q3} Theorem 1.3. In this section we describe his result 
and we introduce some notation.

\begin{dfn} \label{defpath}

A \emph{path} $\gra$ in our graph is a sequence $h^{(m)}\dots h^{(1)}$ such that
$h^{(i)} \in H$ and $h^{(i)}_1 = h^{(i+1)}_0$ for $i =1,\dots,m-1$. 
We define also $\gra_0=h^{(1)}_0$, $\gra_1= h^{(m)}_1$ and we say that 
the degree of $\gra$ is $m$. If $\gra_0=\gra_1$ we say that $\gra$ is 
a closed $B$-path.
The product of paths is defined in the obvious way.

An \emph{admissible path} 
$[\grb]$ in our graph is a sequence  $[ i_{m+1}^{r_{m+1}}
\gra^{(m)}i_{m}^{r_m}\dots\gra^{(1)}i_1^{r_1} ]$, that we write 
between square brackets
such that $i_j \in I$, $\gra^{(j)}$ are $B$-path, $r_j \in \mN$ and 
$\gra^{(j)}_0 = i_j $ and $\gra^{(j)}_1= i_{j+1} $ for $j=1,\dots,m$. 
If $[\grb]$ is as above we 
define $[\grb]_0 = i_1$, $[\grb]_1 = i_{m+1}$.
and the degree of $[\grb]$ is 
$2+\sum_{j=1}^{m+1} r_j +\sum_{j=1}^{m} degree(\gra^{j})$ 
and the product of paths is defined by:
$$
[\grb]\cdot[\grb^{\prime}]=
\begin{cases}
0     &\text{if } [\grb^{\prime}]_1 \neq [\grb]_0 \\
{[\grb i \grb^{\prime}]}  &\text{if } [\grb^{\prime}]_1 = [\grb]_0 =i
\end{cases}
$$

Given a $B$-path $\gra= h^{(m)}\dots h^{(1)}$
and an admissible path $\grb = [ i_{m+1}^{r_{m+1}}
\gra^{(m)}\dots$  $\dots i_1^{r_1} ]$
we define an evaluation of $\gra$ and $\grb$ on $S$
in the following way: if $s=((B)_{h\in H},\grg,\grd) \in S $ then
\begin{align*}
\gra(s) &= B_{h^{(m)}} \circ \dots \circ B_{h^{(1)}} \in 
		\Hom (V_{\gra_0},V_{\gra_1}), \\
\grb(s) &= \grd_{i_{m+1}} \circ (\grg_{i_{m+1}}\circ \grd_{i_{m+1}})^{r_{m+1}}
	\circ \gra^{(m)}(s) \circ (\grg_{i_{m}}\circ \grd_{i_{m}})^{r_{m}}
 	\circ \dots \circ \\
	&\quad \circ \dots \circ \gra^{(1)}(s) \circ 
	(\grg_{i_{1}}\circ \grd_{i_{1}})^{r_{1}} \circ \grg_{i_1}\in 
		\Hom (D_{\grb_0},D_{\grb_1}).
\end{align*}

According to the last definition, 
to simplify the notation we will use the words built with 
the letters $A_i$, $B_i$, $\grg_i$, $\grd_i$
also to denote the paths:
instead of writing $\gra = h^{(m)}\dots h^{(1)}$
we will write  $\gra = B_{h^{(m)}}  \cdots  B_{h^{(1)}}$ and instead of 
$\grb = [ i_{m+1}^{r_{m+1}} \dots i_1^{r_1} ]$ we will write 
$\grd_{i_{m+1}} (\grg_{i_{m+1}} \grd_{i_{m+1}})^{r_{m+1}}$ $\dots \grg_{i_1}$.
In particular $\grd _{l \ra j} \grg_{i \ra l}$ have to be thought 
also as admissible paths.

The algebra $\wt {\calR}$ is  the vector 
space spanned by the admissible path with the product induced by the product 
of path described above. 
Finally we 
define the associative \emph{algebra $\calR=\wt{\calR} / \calI$ of 
admissible polynomials} where 
$\calI$ is the bisided ideal generated by the elements
$[\gra \theta_i \gra^{\prime}] $ where 
$\gra, \gra^{\prime}$ are $B$-path such that $\gra_0 = i= \gra^{\prime}_1$ 
and 
$$
\theta_i =
\begin{cases}
{[}1] -   [\overline{ar}_1  ar_1 ] &\mif i=1, \\
{[}i] +   [ar_{i-1} \overline{ar}_{i-1}] - [\overline{ar}_i  ar_i ] 
&\mif i = 2, \dots, n-2 ,\\
{[}n-1] + [ar_{n-2} \overline{ar}_{n-2}] &\mif i=n-1 .
\end{cases}
$$
If $f$ is an element of $\calR$ 
and there exist $i,j \in I$ such that: 
$$
f = \sum _{\grb \tc \grb_0=i \quad \grb_1=j} a_{\grb} [\grb]
$$
we say that $f$ is of type $(i,j)$.
\end{dfn}
\begin{oss}
We observe that the evaluation on $S$ is a morphism of algebra from $\wt{\calR}$
to the algebra defined by the morphisms of the category of vector spaces.
We observe also that the evaluation of $\calR$ on elements of $\grL$
is well defined. Moreover if $f$ is of type $(i,j)$ we observe that
$f(s) \in \Hom (D_i,D_j)$.
\end{oss}

\begin{teo}[Lusztig, \cite{Lu:Q3} Theorem 1.3]\label{generatoriinvarianti}
The ring $\mC[S]^G$ is generated 
by the polynomials:
\begin{alignat*}{2}
 s& \longmapsto \Tr\left( \gra(s) \right) &\quad&\mfor \gra 
	\text{ a closed $B$-path}, \\
 s& \longmapsto \grf\left(\grb(s)\right) &\quad&\mfor \grb 
	\text{ an admissible path and } 
	\grf \in \left(\Hom(D_{\grb_0},D_{\grb_1}) \right)^* .
\end{alignat*}
\end{teo}

In the case of graph of type $\mathsf{A}_{n-1}$ 
it is easy to see that the first type of 
generators are not necessary. The following lemma is easy.

\begin{lem} \label{polinomi}
1) The algebra of admissible path is generated by the following set 
\begin{equation} \label{insiemeP}
{\calP} = \{ \grd _{l \ra j} \grg_{i \ra l} \tc i,j \in \{1,\dots,n-1\} 
\mand l\leq \min (i,j) \}.
\end{equation}

2) $\mC[\grL]^G$  is generated by the polynomials:
$$s \longmapsto \grf(\grb(s)) \qquad \mfor \grb \in \calP \mand \grf \in 
\left( \Hom(D_{\grb_0},D_{\grb_1}) \right)^*.$$

3)
If $(A,B,\grg,\grd) \in \grL$ then it is 
an element of $ \grL^+$ iff for all $ 1 \leq i \leq n-1$ 
$$\Im A_{i-1} + \sum_{j=i}^{n-1} \Im \grg_{j\ra i} = V_i .  $$
\end{lem}

We will need also a criterion to understand when $M(d,v)$ is not empty:
this is furnished by Nakajima's construction of the irreducible representation
of $sl_n$: [Na2] \S 10. To state it we recall that we identify $d$ (resp. $v$)
with elements of the weight (resp. root) lattice  (see Notation 1) and we observe 
that if $\grs \in S_n$ then $\grs(v-d) +d $ is in the root lattice.

\begin{lem} \label{quiverAvuoto}
If $\grs \in S_n$ is such that $\grs(d-v)$ is dominant and 
$\vpri =\grs(v-d)+d$ then 
$$M(d,v) \neq \vuoto \iff M(d,v')\neq \vuoto \iff 
\vpri_i \geq 0 \; \mfor i=1,\dots,n .$$
\end{lem}
\begin{proof}
In the case of $v \in \mN^{n-1}_{\geq 0}$ this is a consequence of
Nakajima's theorem: \cite{Na2} \S 10.

Suppose now that there exists $i$ such that $v_i <0$. It is enough to prove 
that there exists $j$ such that  $v_j'<0$. This is a consequence of the 
following well known remark applied in the case $u=\grs(v-d)= d - v'$:
if $d$ and $u$ are dominant and $d-u \geq 0 $ (that 
is $d-u \in \sum \mN_{\geq 0}\gra_i$ ) then $d-\tau u \geq 0$ for all
$\tau \in S_n$.
\end{proof}


\subsection{The Slodowy's variety}
In this section we recall some definition
about the nilpotent variety and the 
partial flag variety.

\begin{dfn} \label{nilpotenti}
    If $N$ is a natural number and $D$ is a vector 
    space of dimension $N$ we define ${\calN} = {\calN}_N$ to be the variety of 
    nilpotent elements in $gl(D)$. 
    Computing the dimensions of the Jordan blocks of an element of ${\calN}$ 
we obtain a 
    partion of $N$, and this give us a parametrization of 
    the orbits $O_{\lambda }$, 
    for $\lambda $ a partion of $N$, of the action of $GL(D)$ on ${\calN}$. 
    If $x \in {\calN}$ and $x,y,h$ is a $sl_2$ triple in $gl(D)$ we define the 
	$transversal \: slice$
    to the orbit of $x$ in ${\calN}$ in the point $x$ as:
    $${\calS} _x = \{u \in {\calN} \text{ such that } [u-x,y]=0 \}.$$
    Here and in the sequel, using a nonstandard convention, we admit $0,0,0$ 
as an $sl_2$ triple, so that 
    in the case of $x=0$ we have ${\calS} _0 =  {\calN}$.
\end{dfn}

\begin{dfn} \label{definizioneflag}
   For $N$ an natural number , $a = (a_1,\dots,a_n)$ a 
   vector of nonnegative integers 
   such that $a_1+\cdots +a_n=N$, $D$ a vector space of dimension $N$ we define a  
   \emph{partial  flag} of 
   type $a$ of $D$ to be an increasing 
   sequence $F: \{0\}=F_0 \subset F_1 \subset \cdots \subset F_n=D$ of 
   subspaces of $D$ such that 
   $\dim F_i - \dim F_{i-1} = a_i$. We 
  define ${\calF} _a $ to be the $GL(D)$-homogenous variety of 
   partial flags of type $a$. We define also 
\begin{align*}
 {\widetilde {\calN} _a} &= T^{*}{\calF} _a  
\cong \{ (u,F) \in gl(D) \times {\calF} _a 
 \text{ such that } u(F_i) \subset F_{i-1} \}, \\
 \mu_a &: {\widetilde {\calN} _a} \longrightarrow {\calN} 
 \text{ the projection onto the first factor, and} \\
 {\calF} _a ^{x} &= \mu^{-1}_a (x) \text{ for } x \in {\calN}.
\end{align*}
For $N$, $a$, $D$ as above let 
$\gra =( \gra_1 \geq \gra_2 \geq \cdots \geq \gra_n )$ be a permutation of $a$ 
and define the partition
$\lambda_ a= 1^{\gra_1 -\gra_2} 2^{\gra_2 -\gra_3} \cdots n^{\gra_n}$. 
$\lambda_ a$ is a partition of $N$ 
and it is known that if $(u,F)$ is in 
${\widetilde {\calN} _a}$ then $u $ is in the closure of 
$O_{\lambda _a}$. Moreover the map
$$\mu_a : {\widetilde {\calN} _a} \longrightarrow {\overline O_{\lambda _a}}$$ 
is a resolution of singularity and it is 
an isomorphism over $O_{\lambda _a}$. We define
$$
 \slice = {\calS} _x \cap {\overline O_{\lambda _a}}, \qquad \quad 
 \sopraslice = \mu_a ^{-1} ({\calS} _{a,x} ).
$$
We call $\sopraslice$ the \emph{Slodowy's variety}.
As in the case of quiver varieties will turn to be useful to define
$\sopraslice$ also in the case of a partition $a$ with negative elements: 
$\sopraslice = \vuoto$ if there exists $i$ such that $a_i <0$.
\end{dfn}

The following proposition is well known.

\begin{prp}\label{slicevuoto}
Let $x\in \calN _N$ of type $1^{d_1} 2^{d_2}\dots (n-1)^{d_{n-1}}$ and
$a=(a_1,\dots,a_n)$ a partition of $N$ then:

1) $\sopraslice \neq \vuoto \iff x \in {\overline O_{\lambda _a}} \iff
\forall \,1\leq k\leq n \mand 
\forall\, 1\leq i_1 < i_2 < \dots < i_k \leq n$ 
the following inequality holds:
\begin{equation}\label{disslicevuoto}
d_1 +2 d_2+\dots +kd_k + \dots + kd_{n-1} \geq a_{i_1}+\dots + a_{i_k}
\end{equation}

2) If $\sopraslice \neq \vuoto$ then it is a smooth variety of dimension
$\dim Z_{gl}(x) - \dim Z_{gl}(u_a)$, where $u_a$ is an element of $O_{\grl_a}$.
\end{prp}


\subsection{Nakajima's conjecture} \label{parcongettura}
If $d=(d_1,\dots,d_n)$ and $ v=(v_1,\dots,v_n)$ are two $n-1$-tuples of
integers we define the $n$-tuple $a=a(d,v)=(a_1,\dots,a_n)$ by:
\begin{gather*}
	a_1 = d_1+\dots+d_{n-1}-v_1,  \qquad a_{n}=v_{n-1}, \\
 	\mand a_i=d_i+\dots+d_{n-1}-v_i+v_{i-1} \;\mfor i=2,\dots,n-1.
\end{gather*}
We observe that $\sum_{i=1}^{n} a_i = N = \sum_{i=1}^{n-1} id_i$. Moreover
we observe that once $d$ is fixed the map $a$ gives a bijection between 
$n-1$-tuples of integers $v$ and $n$-tuples of integers $a$ such that 
$\sum a_i =N$.  Indeed we have that 
$$v_{n-1} = a_n \quad 
v_i = a_n + \dots + a_{i+1} - d_{i+1} - 2 d_{i+2} \dots - (n-i-1) d_{n-1}
$$
for $i=1,\dots, n-2$.
Now we can state the main proposition of this paper.  We recall that 
$M(d,v)=M^1(d,v)=\vuoto$ if $v_i < 0 $ for some $i$ and 
$\sopraslice = \slice =\vuoto$ if $a_i <0 $ for some $i$.
The theorem was conjectured by Nakajima in  \cite{Na1}.

\begin{teo} \label{congettura}
  Let $v$, $d$, $N$, $a=a(d,v)$ as above. 
  Let $x \in {\calN}$ be a nilpotent element of type 
  $1^{d_1}\cdot\cdots\cdot (n-1)^{d_{n-1}}$
  then there exist isomorphisms of algebraic varieties $\wt{\varphi}$ between 
  $M(d,v)$ and 
  $\sopraslice$, and
  $\varphi_1$  between $M^1(d,v) $ and $\slice$  
  such that $0 \in M^1(d,v)$ goes to $x \in {\calS} _{d,x}$ and  
  the following diagram commutes:
 \begin{equation} \label{diagramma}
 \begin{CD}
 M(d,v)      @>\wt{\varphi}>>     {\widetilde {\calS} _{d,x}  }  \\
 @V{\pi}VV                        @V{\mu_d}VV                    \\
 M^1(d,v)    @>\varphi_1>>        {\calS} _{d,x}
 \end{CD}
\end{equation}
\end{teo}

\begin{oss}
If $M(d,v) \neq \vuoto$ then it is easy to see that 
 $0 \in M^1(d,v)$: anyway this will be a consequence of the proof.
\end{oss}

We begin the proof of the theorem with some remarks on the degenerate 
cases and on the dimension of the varieties $M(d,v)$ and $\sopraslice$.

\begin{lem} \label{dimensioni}
Let $d,v,N,a$ as above and let $\grs \in S_n$ such that 
$\grs(d-v)$ is dominant and $\vpri = \grs(v-d)+d$ then:


1) If there exists  $i$ such that  $a_i < 0$ then $M(d,v) = \vuoto$. 

2) If there exists  $i$ such that  $v_i < 0$ then $\sopraslice = \vuoto$. 

3) If $\sopraslice \neq \vuoto$ then $M(d,v) \neq \vuoto$ and they are two 
smooth varieties of the same dimension.
\end{lem}

\begin{proof} 

1) Is an easy consequence of lemma \ref{polinomi} point 3).

2) If $v_i < 0$ then we have 
\begin{gather*}
N-(a_1+\dots +a_i)= a_n + \dots + a_{i+1} <  \\
< d_{i+1} + 2 d_{i+2} \dots + (n-i-1) d_{n-1} =
N - (d_1 +\dots + i d_i+\dots +i d_n).
\end{gather*}
So $a_1+\dots + a_i > d_1 +\dots + i d_i+\dots +i d_n$ and
$\sopraslice$ is empty by lemma \ref{slicevuoto}.

3) We observe that the Weyl group $S_n$ acts by permutation 
on the $n$-tuple $a$ and that :
\begin{enumerate}
	\item ${\wt{\calS}}_{\grs (a), x}\neq \vuoto \iff \sopraslice\neq \vuoto$, 
	\item $a(d,\grs(v-d)+d = \grs( a(d,v))$.
\end{enumerate}
The first property is clear from proposition \ref{slicevuoto} (indeed with
a little more effort can be checked that ${\wt{\calS}}_{\grs (a), x}\isocan 
\sopraslice$ but we don't need this result).
The second property is a computation that can easily be checked 
for $\grs = (i,i+1)$. 
So by Lemma 8 it is enough to prove that $\sopraslice \neq \vuoto \then 
M(d,v) \neq \vuoto$ when $d-v$ is dominant. 
 If we set $i_1 = 1, \dots, i_k=k$ in  
the inequality \eqref{disslicevuoto} we obtain $v_k \geq 0$ for $k=1,\dots ,n-1$ 
and by  lemma \ref{quiverAvuoto} $M(d,v) \neq \vuoto$.

 Finally by a result of Nakajima \cite{Na1}, \cite{Na2} if $M(d,v)$ is not empty 
then it is a smooth 
variety of dimension $2 \trasp{v} w - \trasp{v} C v$
and the equality of dimensions is an easy consequence of \ref{slicevuoto}.
\end{proof}


\section{Definition of the map}


In this section we will define the maps $\varphi_1$ and $\wt{\varphi}$ in the case 
$v_i , a_i \geq 0 $ for each $i$.

\begin{lem}[Nakajima, \cite{Na1}] \label{caso semplice}
If $N \geq v_1 \geq \cdots \geq v_{n-1}$ and if $d=(N,0,\dots,0)$ then the \
conjecture is true. 
In this case we have $D=D_1$ and 
$M(d,v) \simeq \wt{{\calN} _a }$ and \
$M^1(d,v) \simeq {\overline O_{\lambda _a}}$ and $M^0(d,v)$ is the closure
of a nilpotent orbit.
\end{lem}
\begin{proof}
The proof is given in \cite{Na1}, but in 
that case Nakajima consider the inverse condition of stability 
so we remind the definition of 
the isomorphism in our case and we  give a sketch of the proof. 
The isomorphism between $M(d,v)$ and ${\calF} _a $ is given by:
$$(A,B,\grg,\grd) \longmapsto ( \grd_1 \grg_1, \{0\} 
\subset \ker \grg_1 \subset \ker A_1 \grg_1 \subset \cdots \subset 
\ker A_{n-1}\cdots A_1 \grg_1)$$
The map between  $M^1(d,v)$ and $ {\overline O_{\lambda _a}}$ or between 
$M^0(d,v)$ and $ \calN$ is given by
$(A,B,\grg,\grd) \mapsto  \grd_1 \grg_1$.
Once the map is defined on $M(d,v)$ it is easy to check that it is bijective 
and that it is
$GL(D)$ equivariant. Now we know that
the map $\mu_a$ is a resolution of singularity and that it is an isomorphism 
over $O_{\lambda_a}$ 
which is a homogenous space. Now by bijectivity and 
equivariance we see that the map 
we have defined must be an isomorphism over this set. 
Now we can prove the lemma by Zarisky main theorem and the 
normality of the closures of 
nilpotent orbits proved by Kraft and Procesi \cite{KraPro}.
\end{proof}

To treat the general case we use the lemma above in the following way. 
Let $d,v,a,\lambda_a$ be given 
as in theorem \ref{congettura} and define 
$\wt{d} _i= 0$ if $i >1$ and $\wt{d}_1= N=\sum_{j=1}^{n-1} jd_j$, 
$\wt{v} _i= v_i + \sum_{j=i+1}^{n-1} (j-i) d_j $. 
We observe that by the lemma above 
$M(\wt{v}, \wt{d}) = T^* {\calF} _a $ and 
$M^1(\wt{v}, \wt{d}) = {\overline O_{\lambda _a}}$. So we 
can think $\slice$ and $\sopraslice$ as subvarieties of 
$M^1(\wt{v}, \wt{d})$ and $M(\wt{v}, \wt{d})$:
\begin{align*}
\slice &=    
p_0 \Bigl( \bigl\{ 
(\wt{A}, \wt{B}, \wt{\grg}, \wt{\grd}) \in \grL(\wt{v}, \wt{d})
\tc [\wt{\grd}_1 \wt{ \grg}_1 -x,y]=0 \bigr\} \Bigr) \cap M^1(\wt v ,\wt d ), \\
\sopraslice &= p \Bigl( \bigl\{ (\wt{A}, \wt{B}, \wt{\grg}, \wt{\grd}) \in 
\grL^+(\wt{v}, \wt{d}) \: | \: [\wt{\grd}_1 \wt{ \grg}_1 -x,y]=0 \bigr\}  \Bigr).
\end{align*}
So we can construct our map by giving a map from $\grL(d,v)$ to $\grL(\wt{v}, \wt{d})$. 
Let us begin with the definition of  $\wt{V}$ and $\wt{D}$. Let  
$D_i^{(j)}$ be an isomorphic copy of $D_i$. We define:
\begin{subequations} \label{eq:tilde}
\begin{align}
   \wt{D_1} &= \wt{D} =  \bigoplus_{1\leq k \leq j \leq n-1}  D_j^{(k)}   \\
   \wt{V_i} &=  V_i \oplus \bigoplus_{1 \leq k \leq j-i  \leq n-i-1}  D_j^{(k)}
\end{align}
\end{subequations}
We will use also the following conventions:  $\wt{V_0} = \wt{D_1}$,  
$ \wt{A}_0 =\wt{\grg}_1$, $\wt{B}_0 = \wt{\grd}_1$ 
and we define the following subspaces
of $\wt V _i$:
$$                                           
     D _i ^{\prime} = \bigoplus_{\substack{  
                                     i+1 \leq j \leq n-1 \\
                                     1 \leq k \leq   j-i}}  D_j^{(k)}  
     \qquad
     D _i ^{+} = \bigoplus_{\substack{
                                     i+2 \leq j \leq n-1 \\
                                     2 \leq k \leq   j-i}} D_j^{(k)}   
     \qquad
     D _i ^{-} = \bigoplus_{\substack{ 
                                     i+2 \leq j \leq n-1 \\
                                     1 \leq k \leq   j-i-1}}  D_j^{(k)}      
$$
We consider the group $GL(V)$ as the subgroup of $GL(\wt V)$ acting as the identity map on $D_i^{\prime} $ 
and mapping $V_i$ into $V_i$.  
We will always think at the maps $\wt A_i, \wt B_i$ as a block-matrix  with respect to
the given decomposition of $\wt V $, $\wt D$ and when we use a projection on
one of our subspaces, it will be a projection with respect  to the given 
decompositions \eqref{eq:tilde}.
To simplify the notation we give a name to these blocks:
\begin{equation} \label{blocchi}
\begin{alignedat}{2}
 \pi_{D^{(h)}_j} \wt A_i|_{D^{(\hp)}_{\jp}}    &=  t^{\jp,\hp}_{i,j,h}     
 &    \qquad   \pi_{D^{(h)}_j} \wt B_i|_{D^{(\hp)}_{\jp}}    
 &=  s^{\jp,\hp}_{i,j,h}  \\ 
  \pi_{D^{(h)}_j} \wt A_i|_{V_i }               &=  t^{V}_{i,j,h}           &    \qquad   \pi_{D^{(h)}_j} \wt B_i|_{V_{i+1} }           &=  s^{V}_{i,j,h}        \\
  \pi_{V_{i+1}}   \wt A_i|_{D^{(\hp)}_{\jp}}    &=  t^{\jp,\hp }_{i,V}      &    \qquad   \pi_{V_{i}}     \wt B_i|_{D^{(\hp)}_{\jp}}    &=  s^{\jp,\hp}_{i,V}    \\  
  \pi_{V_{i+1}}   \wt A_i|_{V_i }               &=  a_i                     &    \qquad   \pi_{V_{i}}     \wt B_i|_{V_{i+1} }           &=  b_i   
\end{alignedat}
\end{equation}

We define also $(x_i,y_i,[x_i,y_i])$ to be the following special $sl_2$ triple 
of $sl(D_i ^{\prime})$:
\begin{align*}
x_i\bigr|_{D_j^{(1)}} &= 0   ,                         \\
x_i\bigr|_{D_j^{(h)}} &= Id_{D_j}: D_j^{(h)} \ra D_j^{(h-1)}, \\
y_i\bigr|_{D_j^{(j-i)}} &= 0,                              \\
y_i\bigr|_{D_j^{(h)}} &= h(j-i -h )Id_{D_j}: D_j^{(h)} \ra D_j^{(h+1)} ,
\end{align*}
and we observe that $x=x_0$, 
$y=y_0$, $[x,y]$ is an $sl_2$ triple in $sl(\wt{D})$ of the type required
in the theorem.

We  introduce now a subset of $\grL(\wt v , \wt d)$. 
To do it we give a formal degree 
to the block of our matrices. Indeed we define two different kind of degrees, $deg$ and $grad$, in
the following way:
\begin{align*}
 \degree(t^{\jp,\hp}_{i,j,h} )  &=  \min(h-\hp+1, h-\hp +1 + \jp - j), \\
 \grad(t^{\jp,\hp}_{i,j,h})  &=  2h - 2\hp + 2 + \jp - j  ,         \\
 \degree(s^{\jp,\hp}_{i,j,h} )  &=  \min(h-\hp, h-\hp  + \jp - j) ,    \\
 \grad(s^{\jp,\hp}_{i,j,h})  &=  2h - 2\hp + \jp - j .
\end{align*}

\begin{dfn}
An element $(\wt{A}, \wt{B}, \wt{\grg}, \wt{\grd})$ of $ \grL(\wt v , \wt d)$ 
is called $transversal$ if it satisfies 
the following relations for $ 0 \leq i \leq n-2$:
\begin{equation} \label{relazioniT}
\begin{aligned}
   t^{\jp,\hp}_{i,j,h}  &= 0         & &\mbox{ if } \degree(t^{\jp,\hp}_{i,j,h}) < 0  \\
   t^{\jp,\hp}_{i,j,h}  &= 0         & &\mbox{ if } \degree(t^{\jp,\hp}_{i,j,h}) = 0  \mbox{ and } (\jp,\hp) \neq (j, h+1) \\
   t^{\jp,\hp}_{i,j,h}  &= Id_{D_j}  & &\mbox{ if } \degree(t^{\jp,\hp}_{i,j,h}) = 0  \mbox{ and } (\jp,\hp) = (j, h+1) \\
   t^{V}_{i,j,h}        &= 0         & &   \\
   t^{\jp,\hp}_{i,V}    &= 0         & &\mbox{ if } \hp \neq 1 \\
   s^{\jp,\hp}_{i,j,h}  &= 0         & &\mbox{ if } \degree(s^{\jp,\hp}_{i,j,h}) < 0  \\
   s^{\jp,\hp}_{i,j,h}  &= 0         & &\mbox{ if } \degree(s^{\jp,\hp}_{i,j,h}) = 0  \mbox{ and } (\jp,\hp) \neq (j, h) \\
   s^{\jp,\hp}_{i,j,h}  &= Id_{D_j}  & &\mbox{ if } \degree(s^{\jp,\hp}_{i,j,h}) = 0  \mbox{ and } (\jp,\hp) = (j, h) \\
   s^{V}_{i,j,h}        &= 0         & &\mbox{ if } h \neq j-i  \\
   s^{\jp,\hp}_{i,V}    &= 0         & & 
\end{aligned}
\end{equation}

\noindent and finally if for each  $0 \leq i \leq n-2$:
$$ [ \pi_{D_i^{\prime}} \wt B_i \wt A_i\bigr|_{D_i^{\prime}} - x_i , y_i]=0. $$
We call $T$ the set of transversal data and we call $T^+$ 
the set of stable data which are also transversal.
\end{dfn}

We observe that  $p(T^+) \subset \sopraslice$ and 
$p_0(T) \cap M^1(\wt v, \wt d) \subset \slice$ and we observe also that 
$T$ and $T^+$
are $GL(V)$ invariant closed subset of $\grL$ and $\grL^+$ respectively. 

We will define our maps $\wt {\varphi}$, $\varphi$ by giving a $GL(V)$ 
equivariant map 
$\Phi$ from 
$\grL(d,v)$ to $T$. If $(A,B,\grg,\grd) \in \grL(d,v)$ its image under 
$\Phi$ is an 
element $ (\wt{A}, \wt{B}, \wt{\grg}, \wt{\grd})$ of $T$ such that:
\begin{subequations} \label{phi}
\begin{align}
    a_i &= A_i           &   b_i &= B_i      \label{phi1} \\
    t^{\jp,1}_{i,V}      &=  \grg_{\jp \ra i+1} &     
    s^{V}_{i,j,j-i}      &=  \grd_{i+1 \ra j}   \label{phi2} \\
    t^{\jp,\hp }_{i,j,h} &=  T^{\jp,\hp }_{i,j,h}(A,B,\grg,\grd)&
    s^{\jp,\hp }_{i,j,h} &=  S^{\jp,\hp}_{i,j,h} (A,B,\grg,\grd)  \label{phi3}   
\end{align}
\end{subequations}
where $S^{\jp,\hp}_{i,j,h}$ and $T^{\jp,\hp }_{i,j,h}$ are admissible polynomials 
of type $(j',j)$ (see definition \ref{defpath}).
\begin{oss}\label{oss:untilde}
The conditions  $t^{\jp,1}_{i,V} = \grg_{\jp \ra i+1}$ for $j' > i+1$ and 
$s^{V}_{i,j,j-i}= \grd_{i+1 \ra j}$  for $j > i+1$ are redundant.
Indeed it is easy to see that if 
$(\wt{A}, \wt{B}, \wt{\grg}, \wt{\grd}) \in T$
and $a_i = A_i$, $b_i = B_i$ and 
$t^{i+1,1}_{i,V}=\grg_{i+1}$, $s^V_{i,i+1,i} =\grd_{i+1}$ then  
\eqref{phi2} is satisfied.
We do not give the details of this simple fact because the argument
is completely similar (but much more simple) to the proof of the next lemma.
\end{oss}
\begin{lem} \label{Lemma}
Given $(A,B,\grg,\grd) \in \grL$ there exists a unique 
$(\wt A, \wt B , \wt {\grg}, \wt {\grd}) \in T$ such that 
\eqref{phi1} and \eqref{phi2} are satisfied.
Moreover there exist homogeneous 
admissible polynomials 
$T^{\jp,\hp}_{i,j,h}$ and $S^{\jp,\hp}_{i,j,h}$ 
of type $(j,j')$ such that \eqref{phi3} is satisfied and 
for $\degree >0$ they are homogeneous  polynomials 
of degree equal to $\grad$ of the following form:
\begin{align*}
T^{\jp,\hp}_{i,j,h} &= \lambda^{\jp,\hp}_{i,j,h} \grd_{r \ra j} \grg_{\jp \ra r} + Q^{\jp,\hp}_{i,j,h} \\
S^{\jp,\hp}_{i,j,h} &= \mu^{\jp,\hp}_{i,j,h} \grd_{r\ra j} \grg_{\jp \ra r} + R^{\jp,\hp}_{i,j,h}
\end{align*}
where $r=j+ \hp -h$ and $P$ and $Q$ are admissible polynomials 
that can be expressed as a linear combination of products of 
admissible polynomials of 
degree 
strictly 
less than $\grad$  (at least each monomial of P and Q is a product of 
two admissible polynomials of positive degree)
and $\lambda^{\jp,\hp}_{i,j,h}$, $\mu^{\jp,\hp}_{i,j,h}$ are rational numbers.

2) Moreover for $i=0,\dots,n-2$ and $\deg >0$ the following inequalities hold:
\begin{align*}
 \lambda^{\jp,\hp}_{i,j,h} & >0 \\
 \intertext{for $\hp=1,\quad i+2 \leq \jp \leq n-1  
	    \mand 1 \leq h \leq j-i-1 \leq n-i-2 $,}
 \lambda^{\jp,\hp}_{i,j,h} + \mu^{\jp,\hp-1}_{i,j,h} & >0 \\
 \intertext{for $1 < \hp \leq \jp -i -1 \leq n-i-2 
	    \mand 1 \leq h \leq j-i-1 \leq n-i-2$,}
 \mu^{\jp,\hp}_{i,j,h} & >0  
\end{align*}
for $1 \leq \hp \leq \jp -i -1 \leq n-i-2, \quad h=j-i 
\mand  i+1 \leq j \leq n-1 $.
\end{lem}
\begin{proof}
We prove this lemma by decreasing induction on $i$. 
We prove that once $(A,B,\grg,\grd) \in \grL(d,v)$ is fixed 
then there exist a unique element $ (\wt{A}, \wt{B}, \wt{\grg}, \wt{\grd}) \in T$
such that \eqref{phi1} and \eqref{phi2} are satisfied. Moreover
we give an inductive formula for the computation of this element and from
this formula will be clear that there exist admissible polynomials
as claimed in the lemma. 

For $i=n-2$ we have that $\wt A_{n-2}$ and $\wt B_{n-2}$ are already 
completely defined by relations \eqref{phi} and they verify the relation 
$\wt A_{n-2} \wt B_{n-2} =0$. Now we assume to have constructed 
$T_{j,*,*}^{*,*}$ and $S_{j,*,*}^{*,*}$ for $j \geq i+1$ as stated in the lemma 
such that $\wt A_j, \wt B_j$ verify the relations requested
to be in $T$. We prove that there exist unique  
$T_{i,*,*}^{*,*}$ and $S_{i,*,*}^{*,*}$ such that: 
\begin{equation}  \label{equazioni}
[ \pi_{D_i^{\prime}} \wt B_i \wt A_i|_{D_i^{\prime}} - x_i ,y_i]=0
\mand
\wt A_i \wt B_i =\wt B_{i+1} \wt A_{i+1},
\end{equation}
and we prove also that they have the required form.
First we observe that the following equations are satisfied by relations \eqref{relazioniT} and
\eqref{phi}:
\begin{gather}
     \pi|_{V_{i+1}} \wt A_i \wt B_i |_{V_{i+1}}     = A_iB_i+\grg_{i+1}\grd_{i+1} 
                  = B_{i+1}A_{i+1} = \pi|_{V_{i+1}} \wt B_{i+1} \wt A_{i+1} |_{V_{i+1}}        \notag \\
     \pi|_{V_{i+1}} \wt A_i \wt B_i |_{D_{j}^{(h)}} = \grd _{h,1} \grg_{j\ra i+1} 
                 = B_{i+1} \grg_{j\ra i+2}= \pi|_{V_{i+1}} \wt B_{i+1} \wt A_{i+1} |_{D_j^{(h)}} \notag \\
     \pi|_{D_{j}^{(h)}} \wt A_i \wt B_i |_{V_{i+1}} = \grd _{h,1} \grd_{i+1 \ra j} 
                 = \grd_{i+2\ra j} A_{i+1} = \pi|_{D_{j}^{(h)}} \wt B_{i+1} \wt A_{i+1} |_{V_{i+1}} \notag
\end{gather}
Let $L =\wt B_{i+1} \wt A_{i+1}$,
$M =\wt A_i \wt B_i$ and
$N =\pi_{D_i^{\prime}} \wt B_i \wt A_i|_{D_i^{\prime}} - x_i$, and as we have done in 
 \eqref{blocchi} we define the blocks $L^{\jp, \hp}_{j,h}$, $M^{\jp, \hp}_{j,h}$ and $N^{\jp, \hp}_{j,h}$.
So we can give the following formulation to equations \eqref{equazioni}:
\begin{align} 
 M^{\jp, \hp}_{j,h} &= L^{\jp, \hp}_{j,h}    \label{ML}\\
     \intertext{for $1 \leq \hp \leq \jp -i-1 \leq n-i-2 \mand 1 \leq h \leq j -i-1 \leq n-i-2 $, } 
 N^{\jp, \jp-i}_{j,h} &= 0              \label{N1} \\
     \intertext{for  $1+i \leq  \jp \leq n-1 \mand 1 \leq h \leq j -i-1 \leq n-i-2 $, }  
 N^{\jp, \hp}_{j,1}   &= 0           \label{N2} \\
    \intertext{for $1+i \leq   j  \leq n-1 \mand 2 \leq \hp \leq \jp -i \leq n-i-1 $, and } 
 \hp(\jp-i-\hp)N^{\jp, \hp+1}_{j,h+1} &= h(j-i-h) N^{\jp, \hp}_{j,h}    \label{N3}
\end{align}
for $1 \leq \hp \leq \jp -i-1 \leq n-i-2 \mand 1 \leq h \leq j -i-1 \leq n-i-2$. 

Now we give a degree $\deg$ and a degree 
$\grad$ also to these new blocks, in the following way:
$$ \deg (L^{\jp, \hp}_{j,h})= \deg (M^{\jp, \hp}_{j,h}) =\deg (N^{\jp, \hp}_{j,h}) = 
 \min (h-\hp+1,h-\hp+1+\jp-j),
$$
$$ \grad (L^{\jp, \hp}_{j,h})= \grad (M^{\jp, \hp}_{j,h}) =\grad (N^{\jp, \hp}_{j,h}) =
 2h-2\hp+2+\jp-j.
$$
Since 
$\min(m-\hp+1,m-\hp+1+\jp -j) + \min(h-m,h-m+l-j) \leq \min(h-\hp+1,h-\hp+1+\jp-j)$ and
$\min(m-\hp,m-\hp+\jp -j) + \min(h-m+1,h-m+1+l-j) \leq \min(h-\hp+1,h-\hp+1+\jp-j)$ 
we have that $\deg$
and $\grad$ behaves well under composition; that is:
\begin{align*}
\deg ( S^{l,m}_{i+1,j,h}) + \deg (T^{\jp, \hp}_{i+1,l,m})   & \leq  
	\deg (L^{\jp, \hp}_{j,h}) \\
\deg ( T^{l,m}_{i,j,h})   + \deg (S^{\jp, \hp}_{i,l,m})     & \leq  
	\deg (M^{\jp, \hp}_{j,h}) \\
\deg ( S^{l,m}_{i,j,h})   + \deg (T^{\jp, \hp}_{i,l,m})     & \leq  
	\deg (N^{\jp, \hp}_{j,h}) \\
\grad ( S^{l,m}_{i+1,j,h}) + \grad (T^{\jp, \hp}_{i+1,l,m})   & =  
	\grad (L^{\jp, \hp}_{j,h}) \\
\grad ( T^{l,m}_{i,j,h})   + \grad (S^{\jp, \hp}_{i,l,m})     & =  
	\grad (M^{\jp, \hp}_{j,h}) \\
\grad ( S^{l,m}_{i,j,h})   + \grad (T^{\jp, \hp}_{i,l,m})     & =  
	\grad (N^{\jp, \hp}_{j,h}) 
\end{align*}
(We observe that in the $N$-case the degree $0$ 
term $-x_i$ respects these rules).
So if the blocks have deg strictly less than $0$ they vanish identically, and if $\deg=0$ then 
to be different from zero we must have $j=\jp$ and $h=\hp-1$ and it is straight forward that also in this
case all the equations are satisfied.  In this way we see that the equations \eqref{N1} and
\eqref{N2} are always satisfied.

Now we argue by induction on $d=\deg >0$ in the following way: we assume to have constructed 
$T^{\jp,\hp}_{i,j,h}$ and $S^{\jp,\hp}_{i,j,h}$ for the blocks with $\deg < d$ 
such that all the
relations \eqref{ML} and \eqref{N3} for blocks with $\deg < d$ 
are satisfied and we prove that 
$T^{\jp,\hp}_{i,j,h}$ and $S^{\jp,\hp}_{i,j,h}$ for blocks of $\deg = d$ 
are uniquely determined 
by the equations \eqref{ML} and \eqref{N3} for blocks of $\deg = d$. 
So we have the following relations:
\begin{align*} 
                       M^{\jp, \hp}_{j,h} &= L^{\jp, \hp}_{j,h}    \\
     \hp(\jp-i-\hp)N^{\jp, \hp+1}_{j,h+1} &= h(j-i-h) N^{\jp, \hp}_{j,h}    
\end{align*}
for $1 \leq \hp \leq \jp -i-1 \leq n-i-2 \mand 1 \leq h \leq j -i-1 \leq n-i-2$ and 
$\min (h-\hp+1,h-\hp+1+\jp-j) =d>0$.  By induction hypothesis we have under this assumptions
on $j,\jp,h,\hp,d$ the following formulas:
\begin{align*}
   L^{\jp, \hp}_{j,h} &= \nu_h \grd_{r \ra j} \grg_{\jp \ra r} + C^{\jp, \hp}_{j,h}  \\
   M^{\jp, \hp}_{j,h} &= S^{\jp, \hp}_{i,j,h+1} + T^{\jp, \hp}_{i,j,h} + D^{\jp, \hp}_{j,h} \\
   N^{\jp, \hp}_{j,h} &= E^{\jp, \hp}_{j,h}+
               \begin{cases} \notag
                   T^{\jp, \hp}_{i,j,h}  &\text{ if } \hp=1 \\
                   T^{\jp, \hp}_{i,j,h} + S^{\jp, \hp -1}_{i,j,h} &\text{ if }  1< \hp \leq \jp -i -1
               \end{cases} \\
   N^{\jp, \hp+1}_{j,h+1} &= F^{\jp, \hp}_{j,h}+
               \begin{cases} \notag
                   S^{\jp, \hp}_{i,j,h+1}  &\text{ if } h=j-i-1 \\
                   T^{\jp, \hp+1}_{i,j,h+1} + S^{\jp, \hp }_{i,j,h+1} &\text{ if }  1 \leq h < j -i -1
               \end{cases}            
\end{align*}
where $r=j+\hp-h$, $C^{\jp, \hp}_{j,h}$, $D^{\jp, \hp}_{j,h}$, $E^{\jp, \hp}_{j,h}$, $F^{\jp, \hp}_{j,h}$ are 
admissible polynomials 
that by induction we already know and that are a linear combination of  products of 
admissible polynomials of degree 
strictly less than $\grad$, and 
$$ \nu_h = \begin{cases} \notag
              1  &\text{ if } \hp =1 \mand h=j-i-1 \\
              \lambda^{\jp, \hp}_{i+1,j,h}  &\text{ if } \hp =1 \mand h <j-i-1 \\
              \mu^{\jp, \hp-1}_{i+1,j,h}    &\text{ if } h=j-i-1 \mand \hp > 1 \\
              \lambda^{\jp, \hp}_{i+1,j,h} + \mu^{\jp, \hp-1}_{i+1,j,h} &\text{ if } \hp >1 \mand h <j-i-1
           \end{cases}   
$$
In any case by induction hypothesis we see that $\nu_h$ is a positive rational number. We observe that 
also the numbers $\hp(\jp-i-\hp)=\alpha_h$, and $h(j-i-h)=\beta_h$ are positive rational numbers. 
Now we group together all the equations with the same $j$ and the same $\jp$ and we solve them
altogether. Once we have fixed  $j$ and $\jp$ the relations between indeces can be written in this form:
$h_0 \leq h \leq h_1 \mand \hp = k + h$,
where $h_1 = j-i-1 $ and:
$$
h_0 = \begin{cases} \notag
            d           &\text{ if } \jp \geq j \\  
            d + j - \jp &\text{ if } \jp < j
      \end{cases} \quad, \qquad
k =   \begin{cases} \notag
            1-d       &\text{ if } \jp \geq j \\  
            1+\jp-j-d &\text{ if } \jp < j  
      \end{cases}
$$ 
We observe also that once 
$j,\jp,d$ are fixed also $r= j+\hp-h$ is fixed. 
Now to write our sistems of equations in a more readable way we introduce the following variables:
$X_h=  T^{\jp, \hp}_{i,j,h}  $ and 
$Y_h=  S^{\jp, \hp}_{i,j,h+1} $; and we observe that variables involved exhaust all
the unknown blocks of type $T^{\jp, *}_{i,j,*}$ and $S^{\jp, *}_{i,j,*}$ of $\deg =d$ that is what we
we want to construct. 
So we can write the equations \eqref{ML} and \eqref{N3} in the following way:
\begin{equation} \label{sistema1}
\begin{aligned}
X_{h_0}+Y_{h_0} &= \nu _{h_0}  \grd_{r \ra j} \grg_{\jp \ra r} +P_{1,{h_0}} \\
&\dots       \\
X_{h_1}+Y_{h_1} &= \nu _{h_1}  \grd_{r \ra j} \grg_{\jp \ra r} +P_{1,{h_1}}
\end{aligned}
\end{equation}
and
\begin{equation} \label{sistema2}
\begin{aligned}
\alpha_{h_0}   (Y_{h_0}+X_{h_0 +1}) &= \beta _{h_0} X_{h_0} +P_{2,{h_0}} \\
\alpha_{h_0+1} (Y_{h_0+1}+X_{h_0 +2}) &= \beta _{h_0+1} (Y_{h_0}+X_{h_0 +1})  +P_{2,{h_0+1}} \\
&\dots  \\
\alpha_{h_1-1} (Y_{h_1-1}+X_{h_1}) &= \beta _{h_1-1} (Y_{h_1-2}+X_{h_1-1}) +P_{2,{h_1-1}} \\
\alpha_{h_1  } Y_{h_1} &= \beta _{h_1} (Y_{h_1 -1}+X_{h_1})     +P_{2,{h_1}}
\end{aligned}
\end{equation}
where $P_{*,*}$ are known polynomials of degree equal to the $\grad=2\deg + |j-\jp|$ of our blocks and that are 
a linear combination of products of polynomials that have degree strictly less than $\grad$.
This system has a unique solution: first we use the equations \eqref{sistema2} to give an expression of 
$Y_{h_1}+X_{h_1+1}$ in terms of $X_{h_0}$ then we sum all the equations \eqref{sistema1} and we 
obtain a formula for $X_{h_0}$ and then we see that we can determine all the others $X_h$ and $Y_h$. 
We observe also that equations \eqref{sistema1} and \eqref{sistema2} give an inductive formula for
the coefficients $\lambda^{\jp,\hp}_{i,j,h}$ and $\mu^{\jp,\hp}_{i,j,h}$. 
Indeed they are the coefficient of
the term $\grd_{r \ra j} \grg_{\jp \ra r}$ in the polynomials $T^{\jp,\hp}_{i,j,h}$ and $S^{\jp,\hp}_{i,j,h}$ 
above so they solve
the same systems \eqref{sistema1} and \eqref{sistema2} but with the constant coefficients $P_{*,*}$ 
equal to zero. So if we use the same variables $X$ and $Y$ for $\lambda$ and $\mu$, we obtain from system
\eqref{sistema2} the following formulas:
\begin{align*}
Y_{h_0}+X_{h_0 +1} &= \rho_{h_0}  X_{h_0} \\
& \dots  \\
Y_{h_1 -1}+X_{h_1} &= \rho_{h_1-1} X_{h_0} \\
Y_{h_1} &= \rho_{h_1} X_{h_0} 
\end{align*}
where $\rho_h$ are positive rational numbers. We observe that the coefficients of the point 2) 
of the lemma are just, with our convention $X_{h_0}, (Y_{h_0}+X_{h_0 +1}), \dots  , Y_{h_1}$. 
So it is enough to
prove that $X_{h_0} >0$. But summing the equations in system \eqref{sistema1} we obtain:
$$X_{h_0}= \frac {\nu_{h_0} + \cdots + \nu_{h1} } { 1+\rho_{h_0}+\cdots+\rho_{h_1}} $$
which is a positive rational number and the lemma is proved.
\end{proof}
\begin{oss} \label{costruzione}
The lemma above show, how is possible to define the map $\Phi$ from 
$\grL(d,v)$ to $T$. An inverse of $\Phi$ is given in
the following way. Take $(\wt{A}, \wt{B}, \wt{\grg} =\wt A _0, \wt{\grd}= \wt B _0) \in T$ and define
$\Phi^{-1}(\wt{A}, \wt{B}, \wt{\grg}, \wt{\grd})=(A,B,\grg,\grd)$, where $A_i = \pi_{V_{i+1}} \wt A _i|_{V_i}$, 
$B_i = \pi_{V_{i}} \wt B _i|_{V_{i+1}}$, $\grg_i=\pi_{V_{i}} \wt A _{i-1}|_{D_{i}^{(1)}}$ and 
$\grd_i=\pi_{D_{i}^{(1)}} \wt B _{i-1}|_{V_{i}}$. It is clear that the new data 
is in $\grL(d,v)$ 
and it also clear that $\Phi^{-1} \circ \Phi = Id_{\grL(d,v)}$. 
The relation $\Phi \circ \Phi^{-1} = Id_T$
follows from the unicity proved in the lemma.
To be more precise it follows 
from the unicity proved in the lemma 
and remark \ref{oss:untilde}.
\end{oss}
\begin{lem} \label{isoZT}
1) $\Phi : \grL(d,v) \lra T$ is a $GL(V)$-equivariant isomorphism.

2) $\Phi(z) \in T^+ \iff z \in \grL^+(d,v)$ and 
$\Phi|_{\grL^+} : \grL^+(d,v) \lra T^+$ is a 
$GL(V)$-equivariant isomorphism
\end{lem}
\begin{proof}
We have just proved 1). To prove 2) we observe that in the case of $(\wt v, \wt d)$ the 
stability condition is equivalent to $\wt A_i$ is an epimorphism for $i=0,\dots,n-2$. We observe also
that if $(\wt{A}, \wt{B}, \wt{\grg}, \wt{\grd})=\Phi(A,B,\grg,\grd) \in T$ we have that $\wt A _i|_{D_i^+}$
is an isomorphism onto $D_{i+1}^{\prime}$. 
Since 
$V_i \oplus D_{i+1}^{(1)} \oplus
\cdots \oplus D_{n-1}^{(1)}$ 
is a complementary space of $D_i^+$ and that $V_{i+1}$ is a complementary 
space of $D_{i+1}^{\prime}$, we conclude, by \eqref{relazioniT} and
\eqref{phi}, that the stability condition in our case is equivalent to 
$A_i \oplus \grg_{i+1}\oplus \cdots \oplus \grg_{n-1 \ra i+1} :V_i \oplus D_{i+1}^{(1)} \oplus
\cdots \oplus D_{n-1}^{(1)} \lra V_{i+1}$ is an epimorhpism for $i=0,\dots, n-2$; which is 
exactly the condition of lemma \ref{polinomi} point 3) 
for the stability of $(A,B,\grg,\grd)$.
\end{proof}

\begin{dfn}
As observed $\Phi$ is a $GL(V)$-equivariant morphism, so we can define $\varphi_0$ and $\wt{\varphi}$ 
as the maps making the following diagrams commute:
\begin{equation*}
\begin{CD}
 \grL(d,v)      @>\Phi>>             T                    \\
 @V{p_0}VV                        @V{p_0}VV                     \\
 M_0(d,v)    @>\varphi_0>>          M_0(\wt v, \wt d)
 \end{CD}
\qquad  \qquad
 \begin{CD}
 \grL^+(d,v)     @>\Phi>>             T^+                   \\
 @V{p}VV                        @V{p}VV                     \\
 M(d,v)      @>\wt{\varphi}>>     \sopraslice
 \end{CD}
\end{equation*} 
and if we set  $\varphi_1=\varphi_0|_{M^1(d,v)}$ we
observe that by definition the diagram \eqref{diagramma} commutes, and that 
$\Im \varphi_1 \subset \mu_d(\sopraslice)=\slice$.
\end{dfn}

\begin{cor}
Let $a,d,v,N$ as in section \ref{parcongettura} then 
$$M(d,v) = \vuoto \iff \sopraslice = \vuoto$$
\end{cor}
\begin{proof}
After lemma \ref{dimensioni} we have only to prove that 
$M(d,v) \neq \vuoto \then \sopraslice \neq \vuoto$
but this is clear since we have constructed a map from $M(d,v)$ to $\sopraslice$.
\end{proof}
\section{Proof of Theorem 12}
\begin{lem} \label{iniettivita}
Let $(\wt{A}, \wt{B}, \wt{\grg}, \wt{\grd}) \in T$ and $\wt{g} \in GL(\wt{V})$ then
$$ \wt g (\wt{A}, \wt{B}, \wt{\grg}, \wt{\grd}) \in T 
\Longrightarrow \exists  g \in GL(V) \mbox{ such that }
\wt g (\wt{A}, \wt{B}, \wt{\grg}, \wt{\grd}) =  g (\wt{A}, \wt{B}, \wt{\grg}, \wt{\grd})$$
\end{lem}
\begin{proof}
We prove first that $\wt g _i (V_i)=V_i$ and $\wt g _i (D_i^{\prime}) = D_i^{\prime}$.
To prove it we introduce for $i=0,\dots, n-2$,
$l=0,\dots,n-2-i$ and $h=0,\dots,n-2-i-l$ the following subspaces of $\wt V _i$:
$$
 D _i ^{l,(h)} = \bigoplus_{\substack{  
                                     0 \leq \hp \leq h \\
                                     i+1+l+\hp \leq j \leq n-1} }  
					D_j^{(j-i-\hp)}. 
$$                
We prove that $\wt g_i (D _i ^{l,(h)})=D _i ^{l,(h)}$. Indeed we observe 
that if $(\wt{A}, \wt{B}, \wt{\grg}, \wt{\grd}) \in T$ then $\wt A_i|_{D _i ^{l,(h)}}$ is an isomorphism 
onto $D _{i+1} ^{l-1,(h)}$ for $l \geq 1$. So we can argue by induction on $i$, taking as 
first step the trivial case $i=0$, that $\wt g_i (D _i ^{l,(h)})=D _i ^{l,(h)}$. We observe that 
$D _i ^{0,(n-i-2)}=D _i ^{\prime}$ and so the we have proved $\wt g_i (D _i ^{\prime})=D _i ^{\prime}$.
Now we observe that if $(\wt{A}, \wt{B}, \wt{\grg}, \wt{\grd}) \in T$ then 
$\pi _ {D _i ^{-}} \wt B _i|_{D _{i+1} ^{\prime}}$ is an isomorphism and that
$ \wt B _i({V_{i+1}}) \subset D _i ^{0,(0)} \oplus V_i$. Since ${D _i ^{-}} \oplus V_i$ 
is the complementary subspace, respect our decomposition, of $D _i ^{0,(0)}$  and 
$\wt g_i (D _i ^{0,(0)}) = D _i ^{0,(0)}$ we can conclude that $\wt g _{i+1} (V_{i+1})=V_{i+1}$.

Now we consider $g_i = \wt g_i|_{V_i}$ and we prove that 
$\wt g (\wt{A}, \wt{B}, \wt{\grg}, \wt{\grd}) =  g (\wt{A}, \wt{B}, \wt{\grg}, \wt{\grd})$. 
Arguing as in remark \ref{costruzione} we see that it is enough to prove that the 
$a_i$, $b_i$ and $t_{i,V}^{i+1, 1}$ and $s_{i,i,1}^V$ of the two elements of $T$ are equal.
By construction we have already proved the equality of the $a_i$ and $b_i$ block.
To prove the equality for the $t$ and the $s$ block we observe that it's enough to prove that
$\wt g _i|{D_{i+1}^{(1)}} =Id_{D_{i+1}^{(1)}}$. To prove it we observe that 
$\wt A _i |_{D _{i} ^{l,(0)}}$ is the identity map from $D _{i} ^{l,(0)}$ to ${D _{i+1} ^{l-1,(0)}} $.
So arguing by induction as above we conclude that 
$\wt g_i |_{D _{i} ^{l,(0)}}$ is the identity map, and finally we observe that ${D_{i+1}^{(1)}} \subset
D _{i} ^{l,(0)}$.
\end{proof}

\begin{lem}
$\varphi_0$ and $\varphi_1$ are closed immersions.
\end{lem}

\begin{proof}
It is enough to prove that $\varphi_0$ is a closed immersion. We observe that 
$M^0(d,v)$ and $M^0(\wt v, \wt d)$ are affine varieties whose coordinate ring is
described in lemma \ref{polinomi}. We will prove that the associate map ${\varphi_0} ^{\sharp}$
between these rings is surjective by showing that is it possible to obtain the 
polynomials
in ${\calP}(d,v)$ from the admissible polynomials for $(\wt d , \wt v)$ through the map $\varphi_0$.
Let us introduce the following $deg$ on the set ${\calP}(v,d)$:
$$\deg(\grd_{r \ra \jp} \grg_{j \ra r}) = \min (j-r+1, \jp -r+1)$$
and we observe that usual degree is given by $\grad(\grd_{r \ra \jp} \grg_{j \ra r})=2 \deg +
|\jp -j|$. 
We will prove the statement by induction on $d=\deg$. 
If $d \leq 0$ then $r \geq j+1$ or $r \geq \jp +1$ and so there are no polynomials 
in the set 
${\calP}$ in this case, and the statement is proved.
If $d > 0$ we  consider the following blocks of degree $d$:
$$(\wt {\grd_1} \wt{\grg_1} )^{\jp,1}_{j,h} =(\wt B_0 \wt A_0)^{\jp,1}_{j,h}=R+ \begin{cases} \notag
        1 \cdot \grd_{j+1-h \ra j} \grg_{\jp \ra j+1-h}  &\text{ if } j=h \\
        \lambda_{0,j,h}^{\jp,1} \cdot \grd_{j+1-h \ra j} \grg_{\jp \ra j+1-h}  &\text{ if } j>h   
     \end{cases}$$
where by induction and lemma \ref{polinomi} $R$ is a 
linear combination of products of monomials with a smaller $deg$. Since by lemma \ref{Lemma}
the coefficient of $\grd_{j+1-h \ra j} \grg_{\jp \ra j+1-h}$ is different from zero we obtain that for any 
$1 \leq h \leq j$ the 
element of ${\calP}(v,d)$,  $\grd_{j+1-h \ra j} \grg_{\jp \ra j+1-h}$, can be obtained as claimed.
But now we observe that this element has $deg = d$ and that all the elements in ${\calP}$
of deg equal to $d$ can be obtained in this way for a good choice of $h$ between $1$ and $j$.
\end{proof}

{\em Proof of theorem} \ref{congettura}.
By the lemma above and the fact that $\mu_d$ and $\pi$ are projective we see that $\wt{\varphi}$ is proper.
By lemmas \ref{isoZT} and \ref{iniettivita}, since by a result of Nakajima (\cite{Na1} \cite{Na2}) 
all the orbits in $\grL^+(v,d)$ and $\grL^+(\wt d , \wt v)$ are closed we see 
that $\wt{\varphi}$ is also injective. Since by 
lemma \ref{dimensioni} $M(v,d)$ and $\sopraslice$ are smooth 
varieties of the same dimension and $\sopraslice$ is connected
we have proved that it is an isomorphism of holomorphic varieties and by consequence is also an isomorphism of algebraic varieties.
In particular $\wt{\varphi}$ is surjective and $\mu_d$ is also surjective, so also $\varphi_1$ is surjective, 
but since it is a closed immersion of reduced varieties over $\mC $ it must be an isomorphism of 
algebraic varieties. Finally $\varphi_0 (0) =x \in \slice$, so by the previous lemma 
$0 \in M^1(v,d)$ and $\varphi_1 (0) =x $. {\em QED}


\begin{oss}
The map  $\wt{\varphi}$ restricted to $\Gamma(v,d)$ take a 
more explicit and simple 
form. Indeed it is easy to see that in this case $\grd$ vanishes 
so we have that all the polynomials
$T$ and $S$ vanish also, and we have an explicit formula for $\wt{\varphi}$. 
\end{oss}

\begin{oss}
In \cite{Na1} is observed that the conjecture does not generalize to diagrams of 
type $E$ and $D$. 
But it is an interesting and more general 
fact (see for example the stratification of quiver varieties constructed by 
Nakajima \cite{Na1},
\cite{Na2} or the remark above)
that some subvarieties can be described as an another quiver variety.
From this point of vied we want to point out that it is possible to 
give an explicit pairs of injective maps 
$\wt{\psi}$ and  $\psi$ from $M(v,d)$ to $M(\wt d , \wt v)$ and from 
$M_0(v,w)$ to $M_0(\wt d , \wt v)$ respectively such that the diagram 
\eqref{diagramma} commute and $\psi(0)=x$. As we said they have an 
explicit formula and so they look more simple than 
$\wt{\varphi}$ and  $\varphi_1$ but their image is not contained 
in $\sopraslice$ and $\slice$ respectively, they "describe" different 
subvarieties.
\end{oss}


\bibliographystyle{amsplain}

\providecommand{\bysame}{\leavevmode\hbox to3em{\hrulefill}\thinspace}


\medskip    

\address{Dipartimento di Matematica, Universit\`a ``La Sapienza'' di Roma,
        P.le A. Moro 5, 00185 ROMA. }
\email{maffei@@mat.uniroma1.it}

\end{document}